\documentclass[12pt]{amsart}
\usepackage{amsfonts,amssymb ,amsthm,latexsym}

\newtheorem{Thm}{Theorem}[section] 
\newtheorem{Cor}[Thm]{Corollary} 
\newtheorem{Prop}[Thm]{Proposition} \newtheorem{Def}[Thm]{Definition}

\newcommand{\no}{\nonumber}  \newcommand{\bea}{\begin{eqnarray}}
\newcommand{\eea}{\end{eqnarray}}

 \def \R{\rm {\bf R}}      \def \Bbb{\mathbb} 
\def \L{L^{2}(\mathbb R)}
\begin{document}

\title{Analyzing the Weyl-Heisenberg Frame Identity}
\author{Peter G. Casazza  and Mark C. Lammers}
\address{Department of Mathematics \\
University of Missouri-Columbia \\
Columbia, MO 65211 and
Department of Mathematics, University
of South Carolina, Columbia, SC 29208
}
\thanks{The first author was supported
by NSF DMS 970618}

\email{pete@math.missouri.edu; lammers@math.sc.edu}


\maketitle

\begin{abstract}
In 1990, Daubechies proved a fundamental identity for Weyl-Heisenberg
systems which is now called the Weyl-Heisenberg Frame Identity.
WH-Frame Identity:    If $g\in W(L^{\infty},L^{1})$, then for all
continuous, compactly supported functions $f$ we have:
\[\sum_{m,n}|<f,E_{mb}T_{na}g>|^{2} = 
\frac{1}{b}\sum_{k}\int_{\Bbb R}\overline{f(t)}f(t-k/b)\sum_{n}
g(t-na)\overline{g(t-na-k/b)}\ dt.\] It has been folklore that the
identity will not hold universally.  We make a detailed study of the
WH-Frame Identity and show: (1)  The identity does not require any
assumptions on $ab$ (such as the requirement that $ab\le 1$ to have a
frame); (2) As stated above, the identity holds for all $f\in
L^{2}(\Bbb R)$; (3)  The identity holds for all bounded, compactly
supported functions if and only if $g\in L^{2}(\Bbb R)$; (4)  The
identity holds for all compactly supported functions if and only if
$\sum_{n}|g(x-na)|^{2}\le B$ a.e.;  Moreover, in (2)-(4) above, the
series on the right converges unconditionally; (5)  In general, there
are WH-frames and functions $f\in L^{2}(\Bbb R)$ so that the series on
the right does not converge (even symmetrically).  We give necessary
and sufficient conditions for it to converge symmetrically;  (6)
There are WH-frames for which the series on the right always converges
symmetrically to give the WH-Frame Identity, but there are functions
for which the series does not converge and we classify when the series
converges for all functions $f\in \L$;  (7)  There   are WH-frames for
which the series always converges, but it does not converge
unconditionally for some functions, and we classify when we have
unconditional convergence for all functions $f$;  and (8)  We show
that the series converges unconditionally for all $f\in L^{2}(\Bbb R)$
if $g$ satisfies the   CC-condition.
\end{abstract}

\section{Introduction} \label{intro}

In 1990, Daubechies \cite{D} proved a fundamental identity for  
Weyl-Heisenberg
systems , which is now called the {\bf Weyl-Heisenberg Frame  
Identity} (or
{\bf WH-frame Identity} for short).  This identity has been extensively
used in the theory and has gone through some small imporvements  
over time.
It has been part of the folklore that the identity does not hold  
universally.
But, until now, it has been a little mysterious as to exactly when  
and where
can one be sure the identity holds.  In this paper we give a  
detailed analysis
of the WH-frame Identity and answer all the relevant questions  
completely.\\

Casazza, Christensen, and Janssen \cite{CCJ} made a detailed study of
Weyl-Heisenberg frames, translation invariant systems and the Walnut
representation of the frame operator.  We will rely heavily here on these
results and the relevant constructions
from \cite{CCJ} using the Zak transform.\\

The authors would like to thank A.J.E.M. Janssen for interesting  
discussions
concerning the results in this paper.

\section{Preliminaries}\label{1}
\setcounter{equation}{0}

In this section we will give the basic results needed throughout the
paper.  We use $\Bbb N, \Bbb Z, \Bbb R, \Bbb C$ to denote the  
natural numbers,
integers, real numbers and complex numbers, respectively.  A {\bf
scalar} is an element of $\Bbb R$ or $\Bbb C$.  Integration is always
with respect to Lebesgue measure.  $L^{2}(\Bbb R)$ will denote the
complex Hilbert space of square integrable functions mapping   $\Bbb
R$ into $\Bbb C$.  A bounded unconditional basis for a Hilbert space
$H$ is   called a {\bf Riesz basis}.  That is, $(f_{n})$ is a Riesz
basis for $H$ if   and only if there is an orthonormal basis $(e_{n})$
for $H$ and an  operator  $T:H\rightarrow H$ defined by $T(e_{n}) =
f_{n}$, for all $n$.  We call $(f_{n})$ a {\bf Riesz basic sequence}
if it is a Riesz basis for its closed linear span.  \\

In 1952, Duffin and Schaeffer \cite{DS} were working on some deep
problems in non-harmonic Fourier series.  This led them to define

\begin{Def}\label{Frame}  A sequence $(f_{n})_{n\in \Bbb Z}$ of elements
of a Hilbert space $H$ is called a {\bf frame} if there are constants
$A,B>0$ such that
\begin{equation}\label{FrameEqn}
A\|f\|^{2}\le \sum_{n\in \Bbb Z}|<f,f_{n}>|^{2}\le B\|f\|^{2},\ \
\text{for all}\ \ f\in H.
\end{equation}
\end{Def}

The numbers $A,B$ are called the {\bf lower} and {\bf upper frame
bounds} respectively.    The frame is a {\bf tight frame} if
$A=B$ and a {\bf normalized tight frame} if $A=B=1$.  A frame is {\bf
exact} if it ceases to be a frame when any one of its elements is
removed.  It is known that a frame is exact if and only if it is a
Riesz basis.  A non-exact frame is called {\bf over-complete}  in the
sense that at   least one vector can be removed from the frame and the
remaining set of   vectors will still form a frame for $H$ (but
perhaps with different frame bounds).  If $f_{n}\in H$, for all $n\in
\Bbb Z$, we call $(f_{n})_{n\in \Bbb Z}$ a {\bf frame sequence} if it
is a frame for its closed linear span in $H$.\\

We will consider frames from the operator theoretic point of view.  To
formulate this approach, let $(e_{n})$ be an orthonormal basis for an
infinite dimensional Hilbert space $H$ and let $f_{n}\in H$, for all
$n\in \Bbb Z$.  We call the operator $T:H\rightarrow H$ given by
$Te_{n}=f_{n}$ the {\bf preframe operator} associated with $(f_{n})$.
Now, for each $f\in H$ and $n\in \Bbb Z$ we have $<T^{*}f,e_{n}> =
<f,Te_{n}>=<f,f_{n}>$.  Thus
\begin{equation} \label{ET}
T^{*}f = \sum_{n}<f,f_{n}>e_{n},\ \ \text{for all} \ \ f\in H.
\end{equation}

By (\ref{ET})
$$
\|T^{*}f\|^{2} = \sum_{n}|<f,f_{n}>|^{2},\ \ \text{for all}\ \ f\in H.
$$
It follows that the preframe operator is bounded if and only if
$(f_{n})$ has a finite upper frame bound $B$.  Comparing this to
Definition \ref{Frame} we have

\begin{Thm}\label{T1}
Let $H$ be a Hilbert space with an orthonormal basis $(e_{n})$. Also
 let $(f_{n})$ be a sequence of elements of $H$ and let $Te_{n}=f_{n}$
 be the preframe operator.  The following are equivalent:

(1)  $(f_{n})$ is a frame for $H$.

(2)  The operator $T$ is bounded, linear and onto.

(3)  The operator $T^{*}$ is an (possibly into) isomorphism called
the {\bf frame transform}.

Moreover, $(f_{n})$ is a normalized tight frame if and only if the
preframe operator is a quotient map (i.e. a co-isometry).
\end{Thm}

The dimension of the kernel of T is called the {\bf excess} of the
frame.  It follows that $S=TT^{*}$ is an invertible operator on   $H$,
called the {\bf frame operator}.  Moreover, we have
$$
Sf = TT^{*}f = T(\sum_{n}<f,f_{n}>e_{n}) = \sum_{n}<f,f_{n}>Te_{n}
=\sum_{n}<f,f_{n}>f_{n}.
$$

A direct calculation now yields
$$
<Sf,f> = \sum_{n}|<f,f_{n}>|^{2}.
$$
Therefore, the {\bf frame operator is a positive, self-adjoint
invertible operator} on $H$.  Also, the frame inequalities
(\ref{Frame}) yield that $(f_{n})$ is a frame with frame bounds
$A,B>0$ if and only if $A\cdot I\le S\le B\cdot I$.  Hence, $(f_{n})$
is a normalized tight frame if and only if $S=I$.  Also, a direct
calculation yields \bea\label{E1} f = SS^{-1}f &=&
\sum_{n}<S^{-1}f,f_{n}>f_{n} \\ &=& \sum_{n}<f,S^{-1}f_{n}>f_{n}\no \\
&=&\sum_{n}<f,S^{-1/2}f_{n}>S^{-1/2}f_{n}.\no \eea

We call $(<S^{-1}f,f_{n}>)$ the {\bf frame coefficients} for $f$.  One
interpretation of equation (\ref{E1}) is that $(S^{-1/2}f_{n})$ is a
normalized tight frame.

We will work here with a particular class of frames called
Weyl-Heisenberg frames.  To formulate these frames, we first need some
notation.  For a function $f$ on $\Bbb R$ we define the operators:

\[
\begin{array}{lll}
\text{Translation:} & T_{a}f(x) = f(x-a), & a\in \Bbb R \\
\text{Modulation:}  & E_{a}f(x) = e^{2{\pi}iax}f(x), & a\in \Bbb R \\
\end{array}
\]

We also use the symbol $E_{a}$ to denote the {\bf exponential
function} $E_{a}(x) = e^{2{\pi}iax}$.  Each of the operators $T_{a},
E_{a}$ are unitary operators on  $L^{2}(\Bbb R )$

In 1946 Gabor \cite{G} formulated a fundamental approach to signal
decomposition in terms of elementary signals.  This method resulted in
{\bf Gabor frames} or as they are often called today {\bf
Weyl-Heisenberg frames}.

\begin{Def}\label{WHS}
If $a,b\in \Bbb R$ and $g\in L^{2}(\Bbb R )$ we call
$(E_{mb}T_{na}g)_{m,n\in \Bbb Z}$ a {\bf Weyl-Heisenberg system} ({\bf
WH-system} for short) and denote it by $(g,a,b)$.  We call $g$ the
{\bf window function}.
\end{Def}

If the WH-system $(g,a,b)$ forms a frame for $L^{2}(\Bbb R)$, we
call this a {\bf Weyl-Heisenberg frame} ({\bf WH-frame} for short).
The   numbers $a,b$ are the {\bf frame parameters} with $a$ being the
{\bf shift parameter}   and $b$ being the {\bf modulation parameter.}
We will be interested in   when there are finite upper frame bounds
for a WH-system.  We call this class of functions the {\bf preframe
functions} and denote this class by {\bf PF}.  It is easily checked
that

\begin{Prop}\label{P6}
The following are equivalent:

(1) $g\in$ {\bf PF}.

(2)  The operator
$$
Sf = \sum_{n}<f,E_{mb}T_{na}g>E_{mb}T_{na}g,
$$
is a well defined bounded linear operator on $L^{2}(\Bbb R)$.
\end{Prop}

A family $(g,a,b)$ with $g\in $PF is called a {\bf preframe
WH-system}.

It is a simple calculation to check the following (see \cite{CL}):

\begin{Prop}\label{sums}
Let $f,g\in \L$ and $a,b\in \Bbb R$.

(1)  We have
$$
\sum_{k\in \mathbb Z}f(t-k/b)g(t-k/b-na)\in L^{1}[0,1/b].
$$

(2)  If $\sum_{k}|f(t-ka)|^{2}\le B$ then for all $n\in \Bbb Z$,
$$
\sum_{k\in \mathbb Z}f(t-ka)g(t-ka-n/b)\in L^{2}[0,1/b].
$$

(3)  If $\sum_{n}|g(t-na)|^{2}\le B$ then $\sum_{n}|g(t-na)g(t-na-k/b)|
\le B$, for all $k\in \Bbb Z$.
\end{Prop}

We next recall the {\bf Wiener amalgam space} $W(L^{\infty},L^{1})$
which consists of all functions $g$ so that for some $a>0$ we have,
$$
\|g\|_{W,a} = \sum_{n\in \Bbb Z}\|g\cdot
{\chi}_{[an,a(n+1))}\|_{\infty} = \sum_{n\in \Bbb Z}\|T_{na}\cdot
{\chi}_{[0,a)}\|_{\infty} < \infty.
$$
It is easily checked that $W(L^{\infty},L^{1})$ is a Banach space
with the above norm.  Also, if $\|g\|_{W,a}< \infty$, for one $a>0$,
then   this norm is finite for all $a>0$.\\

To simplify some of the results we introduce some notation.  For
any $a,b\in \Bbb R$ and $g\in \L$ we let for all $k\in \Bbb Z$,
$$
G_{k}(t) = \sum_{n\in \Bbb Z}g(t-na)\overline{g(t-na-k/b}.
$$
In particular,
$$
G_{0}(t) = \sum_{n\in \Bbb Z}|g(t-na)|^{2}.
$$

Our main tool will be the proof of the WH-frame Identity due to
Walnut \cite{W}.  He eliminated the need for the Poisson summation
formula used by Daubechies in the original proof and obtained a
more general result.

\begin{Thm}\label{WHFI}({\bf WH-Frame Identity.})
If $g\in W(L^{\infty},L^{1})$  and $f\in L^{2}(\Bbb R)$ is
continuous and compactly supported, then
$$
\sum_{n,m\in \Bbb Z}|<f,E_{mb}T_{na}g>|^{2} =
$$
$$
\frac{1}{b}\sum_{k\in \mathbb Z}\int_{\mathbb  
R}\overline{f(t)}f(t-k/b)G_{k}(t)\ dt = F_{1}(f)+F_{2}(f),
$$
where
$$
F_{1}(f) = b^{-1}\int_{R}|f(t)|^{2}G_{0}(t)\  dt,
$$
and \bea F_{2}(f)  &=& b^{-1}\sum_{k\not=
0}\int_{R}\overline{f(t)}f(t-k/b)G_{k}(t)\ dt \no \\
&=&b^{-1}\sum_{k\ge
1}2\text{Re}\int_{R}\overline{f(t)}f(t-k/b)G_{k}(t)\ dt.  \no \eea
\end{Thm}

{\it Proof.}  We are assuming that $f$ is bounded and compactly  
supported so that all the
summations, integrals and interchanges of these below are  
justified.  We define
$$
H_{n}(t)=\sum_{k}f(t-k/b)\overline{g(t-na-k/b)}.
$$
Now, $H_{n}$ is 1/b-periodic, $H_{n}\in L^{2}[0,1/b]$ and
$$
\int_{R}f\cdot \overline{E_{mb}T_{na}g(t)}dt=  
\int_{R}f(t)\overline{g(t-na)}e^{-2{\pi}imbt}dt =
\int_{0}^{1/b}H_{n}(t)e^{-2{\pi}imbt}dt.
$$
Since $(b^{1/2}E_{mb})_{m\in \Bbb Z}$ is an orthonormal basis for  
$L^{2}[0,1/b]$, the Plancherel
formula yields
$$
\sum_{m}|\int_{0}^{1/b}H_{n}(t)e^{-2{\pi}imbt}dt|^{2} =  
b^{-1}\int_{0}^{1/b}|H_{n}(t)|^{2}dt.
$$
Now we compute
$$
\sum_{n}\sum_{m}|<f,E_{mb}T_{na}g>|^{2} =
\sum_{n}\sum_{m}|\int_{R}f(t)\overline{g(t-na)}e^{-2{\pi}imbt}dt|^{2}
$$
$$
=  
b^{-1}\sum_{n}\int_{0}^{1/b}|\sum_{k}f(t-k/b)\overline{g(t-na-k/b)}|^{2}dt
$$
$$
=  
b^{-1}\sum_{n}\int_{0}^{1/b}\sum_{\ell}\overline{f(t-{\ell}/b)}g(t-na-{\ell}/b)\cdot
\sum_{k}f(t-k/b)\overline{g(t-na-k/b)}dt
$$
$$
=  
b^{-1}\sum_{n}\sum_{\ell}\int_{0}^{1/b}\overline{f(t-{\ell}/b)}g(t-na-{\ell}/b)\cdot
\sum_{k}f(t-k/b)\overline{g(t-na-k/b)}dt
$$
$$
= b^{-1}\sum_{n}\int_{R}\overline{f(t)}g(t-na)\cdot  
\sum_{k}f(t-k/b)\overline{g(t-na-k/b)}dt
$$
$$
= b^{-1}\sum_{k}\int_{R}\overline{f(t)}f(t-k/b)\cdot
\sum_{n}g(t-na)\overline{g(t-na-k/b)}dt
$$
$$
= b^{-1}\int_{R}|f(t)|^{2}\cdot \sum_{n}|g(t-na)|^{2}dt \ \ +
$$
$$
 b^{-1}\sum_{k\not= 0}\int_{R}\overline{f(t)}f(t-k/b)\cdot  
\sum_{n}g(t-na)\overline{g(t-na-k/b)}dt.
$$
This completes the first part of the WH-Frame Identity.  The  
equality in the last line follows
by a simple change of variables.
\qed  \vspace{14pt}

To avoid ``technicalities'' we will say that the {\bf WH-frame Identity
holds for a function f} to mean that the series on the left hand side
sums to be finite and is equal to the right hand side sum which
converges unconditionally.  Later, we will discuss different forms
of convergence for the right hand side of the WH-frame Identity.\\

As a consequence of the WH-frame Identity, Casazza, Christensen and
Janssen \cite{CCJ} showed:

\begin{Prop}\label{PP}
Let $a,b\in R$ with $ab\le 1$ and $g\in L^{2}(R)$ and assume that
$$
\sum_{k\in Z}|G_{k}(t)|^{2} \le B, \ \ \text{a.e.}
$$
Then for all bounded, compactly supported functions $f\in L^{2}(R)$ 
the series
$$
Lf = b^{-1}\sum_{k}(T_{k/b}f)G_{k},
$$
converges unconditionally in norm in $L^{2}(R)$.  Moreover,
$$
<Lf,f> = \sum_{m,n\in Z}|<f,E_{mb}T_{na}g>|^{2}.
$$
Finally, if $g\in$ PF, so that the series
$$
Sf = \sum_{m,n}<f,E_{mb}T_{na}g>E_{mb}T_{na}g,
$$
also converges unconditionally in $L^{2}(R)$, we have that
$Lf = Sf$.
\end{Prop}

We will also make use of the CC-condition from \cite{CC}.

\begin{Thm}[CC-Condition]\label{CC}
If $g\in L^{2}(\Bbb R)$, $a,b\in \Bbb R$ and
\begin{equation}\tag{CC}
\sum_{k\in \Bbb Z}|\sum_{n\in \Bbb Z}g(t-na)\overline{g(t-na-k/b)}|  
= \sum_{k\in \Bbb Z}|G_{k}(t)|
\le B,\ \ \text{a.e.},
\end{equation}
then $g\in $ {\bf PF}.
Moreover, if we also have
\begin{equation}\label{ECC1}
\sum_{k\not= 0}|G_{k}(t)|\le (1-{\epsilon})G_{0}(t)\ \ \text{a.e.},
\end{equation}
for some $0<{\epsilon} < 1$,
then $(g,a,b)$ is a WH-frame.
\end{Thm}

We will need a special representation of the frame
operator for WH-frames due to Walnut \cite{W}.

\begin{Thm} \label{Walnut}
Let $a,b>0$ and $g\in W(L^{\infty},L^{1})$ be given.  For each  
$f\in \L$,
the sum $Sf$ converges and is given by
$$
Sf = \frac{1}{b}\sum_{k\in \mathbb Z}T_{k/b}f\cdot G_{k}.
$$
\end{Thm}

The series in Theorem \ref{Walnut} is called the
{\bf Walnut representation of
the frame operator}.  The precise conditions under which the Walnut
representation converges to the frame operator are quite delicate and
were studied in detail in \cite{CCJ}.

\section{Bounded, Compactly Supported Functions and the WH-Frame  
Identity}\label{1}
\setcounter{equation}{0}

We start with a simple observation.

\begin{Prop}\label{gip}
If $a,b\in \Bbb R$ and  $g\in \L$ is bounded and compactly
supported, then $g\in$PF.
\end{Prop}

{\it Proof.}
First, assume that $g$ is supported on $[0,a]$.  Since $g$ is  
bounded above
and compactly supported,
there is a constant $B$ so that
\begin{equation}\label{1}
\sum_{k\in \mathbb Z}|g(t-k/b)|^{2}\le B.
\end{equation}
We define the preframe operator $L:{\ell}_{2}\otimes {\ell}_{2}
\rightarrow \L$ by
$$
L(\sum_{m,n\in \Bbb Z}a_{mn}e_{mn}) = \sum_{m,n\in \Bbb  
Z}a_{mn}E_{mb}T_{na}g,
$$
where $(e_{mn})$ is the natural orthonormal basis of  
${\ell}_{2}\otimes {\ell}_{2}$.
We need to show that $L$ is a bounded operator.  By our assumption on
the support of $g$, we see that $(T_{na}g)_{n\in \mathbb Z}$ are
disjointly supported functions.  Hence,
\begin{equation} \label{2}
\|L(\sum_{m,n\in \mathbb Z}a_{mn}e_{mn})\|^{2} = \sum_{n\in \mathbb Z}
\|\sum_{m\in \mathbb Z}a_{mn}E_{mb}T_{na}g\|^{2}.
\end{equation}
Applying inequality \ref{1} above at the appropriate step, we have
$$
\|\sum_{m\in \mathbb Z}a_{mn}E_{mb}T_{na}g\|^{2} = \int_{\mathbb R}
|\sum_{m\in \mathbb Z}a_{mn}E_{mb}T_{na}g(t)|^{2}\ dt
$$
$$
= \int_{0}^{1/b}|\sum_{m\in \mathbb Z}a_{mn}E_{mb}|^{2}\sum_{k\in  
\mathbb Z}
|g(t-k/b-na)|^{2}\ dt
$$
$$
\le B\int_{0}^{1/b}|\sum_{m\in \mathbb Z}a_{mn}E_{mb}|^{2}\ dt
$$
$$
= B\sum_{m\in \mathbb Z}|a_{mn}|^{2}.
$$
It follows from equation (3.2),
$$
\|\sum_{m\in \mathbb Z}a_{mn}E_{mb}T_{na}g\|^{2}\le \sum_{mn\in  
\mathbb Z}
|a_{mn}|^{2}.
$$
Hence, $L$ is a bounded operator.

For the general case, we observe that $g$ can be written as a finite sum,
say k,
of translates of functions supported on $[0,a]$ and so the preframe  
function
is also bounded in this case by $k\|L\|$.
\qed \vspace{14pt}

\begin{Cor}
If $g\in \L$, then for every bunded, compactly supported function  
$f$ on $\mathbb R$, we have
$$
\sum_{m,n\in \mathbb Z}|<f,E_{mb}T_{na}g>|^{2}< \infty.
$$
\end{Cor}

{\it Proof.}
By Proposition \ref{gip}, if $f$ is bounded and compactly supported then
$(E_{mb}T_{na}f)_{m,n\in \mathbb Z}$ has a finite upper frame  
bound, say $B$.
Now
$$
\sum_{m,n\in \mathbb Z}|<f,E_{mb}T_{na}g>|^{2} = \sum_{m,n\in \mathbb Z}
|<T_{-na}E_{-mb}f,g>|^{2} =
$$
$$
\sum_{m,n\in \mathbb Z}|e^{-2{\pi}imb(x-na)}<E_{mb}T_{na}f,g>|^{2} = 
\sum_{m,n\in \mathbb Z}|<E_{mb}T_{na}f,g>|^{2}
\le B.
$$
\qed \vspace{14pt}

We now present the main result of this section.

\begin{Thm}\label{BCF}
Let $g$ be a measurable function on $\mathbb R$.  The following are
equivalent:

(1)  $g\in \L$.

(2)  The WH-frame Identity holds for all bounded, compactly supported
functions $f$ on $\mathbb R$.
\end{Thm}

{\it Proof.}
$(1)\Rightarrow (2)$:  We assume that $f$ is supported on $[-N,N]$ and
bounded above by $B$.  For a fixed $n\in \mathbb Z$ we consider the
1/b-periodic function
$$
H_{n}(t) = \sum_{k\in \mathbb Z}f(t-k/b)\overline{g(t-na-k/b)}.
$$
Now, the above sum only has $2N$ non-zero terms for each $t\in  
\mathbb R$.
So we can easily follow Walnut's proof the the WH-frame Identity line for
line interchanging the (now finite) sums and integrals until we arrive
at:
$$
\sum_{m,n\in \mathbb Z}|<f,E_{mb}T_{na}g>|^{2} =
$$
$$
\frac{1}{b}\sum_{n\in \mathbb Z}\int_{\mathbb  
R}\overline{f(t)}g(t-na)\cdot
\sum_{k\in \mathbb Z}f(t-k/b)\overline{g(t-na-k/b)}=
$$
$$
\frac{1}{b}\sum_{n\in \mathbb Z}\int_{-N}^{N}\overline{f(t)}g(t-na)\cdot
\sum_{k\in \mathbb Z}f(t-k/b)\overline{g(t-na-k/b)}
$$
To finish the identity, we just need to justify interchanging the  
infinite
sum over $n$ with the finite sum over $k$.  To justify this we  
observe that:
$$
\sum_{n,k}|f(t)|g(t-na)||f(t-k/b)||f(t-na-k/b)| \le
$$
$$
B^{2}\sum_{k=-N}^{N}\sum_{n\in \mathbb Z}|g(t-na)g(t-na-k/b)|.
$$
By Proposition \ref{sums}, we have that
$$
\sum_{n\in \mathbb Z}|g(t-na)g(t-na-k/b)|\in L^{1}[0,a],
$$
and hence
$$
\sum_{n\in \mathbb Z}|g(t-na)g(t-na-k/b)|\in L^{1}[-N,N].
$$
Therefore, we justify the needed interchange of sums and sums with
integrals by the Lebesgue Dominated Convergence Theorem.

$(2)\Rightarrow (1)$:  We do this by contradiction.  If $g$ is not  
square
integrable on $\mathbb R$, then
$$
\|g\|^{2} = \int_{0}^{a}\sum_{n\in \mathbb Z}|f(t-na)|^{2}\ dt = \infty.
$$
Hence, there is some interval $I$ of length $c$ with $0<c<1/b$ so that
$$
\int_{I}\sum_{n\in \mathbb Z}|g(t-na)|^{2}\ dt=\infty.
$$
If we  let $f={\chi}_{I}$, then
the right hand side of the WH-frame Identity becomes
$$
\int_{\mathbb R}|f(t)|^{2}\sum_{n\in \mathbb Z}|g(t-na)|^{2}\ dt
= \int_{I}\sum_{n\in \mathbb Z}|g(t-na)|^{2}\ dt=\infty.
$$
So the right hand side of the WH-frame identity is not a finite
unconditionally convergent series.  i.e.  The WH-frame Identity fails.   
\qed \vspace{14pt}

\section{Compactly supported functions and the WH-Frame  
Identity}\label{2}
\setcounter{equation}{0}

In this section we will drop the hypotheses that our function $f$ has to 
be bounded and discover necessary and sufficient conditions for the
WH-frame Identity to hold.  The conditions are a little stronger than
those required for bounded, compactly supported functions.

\begin{Thm}\label{CF}
Let $g$ be a measurable function on $\mathbb R$.  The following are
equivalent:

(1)  There is a constant $B>0$ so that
$$
\sum_{n\in \mathbb Z}|g(t-na)|^{2}\le B,\ \ \text{a.e.}
$$

(2)  The WH-frame Identity holds for all compactly supported functions
 $f$ on $\mathbb R$.
\end{Thm}

{\it Proof.}
$(1)\Rightarrow (2)$:  If $f$ is compactly supported, we see immediately 
that the sum over k in the right hand side of the WH-frame identity is
a finite sum.  So let
$f_{\ell}(t) = f(t)$ if $|f(t)|\le \ell$ and zero otherwise.  Now,
by Theorem \ref{BCF}
the WH-frame identity holds for all $f_{\ell}$.  That is, for all  
${\ell}\in
\mathbb Z$ we have
$$
\sum_{m,n\in \mathbb Z}|<f_{\ell},E_{mb}T_{na}g>|^{2} =
$$
$$
\frac{1}{b}\sum_{k\in \mathbb Z}\int_{\mathbb  
R}\overline{f_{\ell}(t)}f_{\ell}(t-k/b)G_{k}(t)\ dt
$$

Now we will finish the proof in three steps.
\vspace{14pt}

{\it Step 1:}  We show that
$$
\sum_{m,n}|<f,E_{mb}T_{na}g>|^{2} =
$$
$$
=  
b^{-1}\sum_{n}\int_{0}^{1/b}|\sum_{k}f(t-k/b)\overline{g(t-na-k/b)}|^{2}dt.
$$

To prove Step 1, we let for a fixed $n\in \mathbb Z$,
$$
H_{n}(t) = \sum_{k\in \mathbb Z}f(t-k/b)\overline{g(t-na-k/b)}.
$$
Since the sum on the right hand side above is finite, we can copy
the first steps of the Walnut proof of the WH-frame Identity to
the point of the identity given in Step 1.
\vspace{14pt}

{\it Step 2:}  We show,
$$
\lim_{\ell \rightarrow \infty}\sum_{n}
\int_{0}^{1/b}|\sum_{k}f_{\ell}(t-k/b)g(t-na-k/b)|^{2}
=
$$
$$
=  
b^{-1}\sum_{n}\int_{0}^{1/b}|\sum_{k}f(t-k/b)\overline{g(t-na-k/b)}|^{2}dt=
$$
$$
\sum_{m,n}|<f,E_{mb}T_{na}g>|^{2}.
$$

For step 2,
choose an $N$ so that for all $t\in [0,1/b]$, $f(t-k/b)=0$ for all
$|k|> N$.  Hence, for all $t\in [0,1/b]$ we have
$$
\sum_{n}|\sum_{k}f_{\ell}(t-k/b)g(t-na-k/b)|^{2} \le
$$
$$
\sum_{k=-N}^{N}|f_{\ell}(t-k/b)|^{2}\sum_{n}\sum_{k=-N}^{N}|g(t-na-k/b)|^{2}
\le
$$
$$
\sum_{k=-N}^{N}|f(t-k/b)|^{2}2NB^{2}\in L^{1}[0,1/b].
$$
So, Step 2 follows by the Lebesgue Dominated convergence Theorem.

The following step will complete the proof.
\vspace{14pt}

{\it Step 3:}
$$
\lim_{\ell \rightarrow \infty}\frac{1}{b}\sum_{k\in \Bbb Z}\int_{\Bbb R}
\overline{f_{\ell}(t)}f_{\ell}(t-k/b) G_{k}(t)\ dt =
$$
$$
\frac{1}{b}\sum_{k\in \Bbb Z}\int_{\Bbb R}
\overline{f(t)}f(t-k/b) G_{k}(t)\ dt
$$

For Step 3, note that support $f_{\ell}
\subset$ support $f_{\ell +1}\subset$ support $f$.  Hence, for $k$ fixed
we have:
$$
|f_{\ell}(t)||T_{k/b}f_{\ell}(t)||G_{k}(t)|\uparrow
|f||T_{k/b}f(t)||G_{k}(t)|.
$$
Also, by assumption $|G_{k}(t)|\le B^{2}$.  Since $f\in L^{2}(\Bbb  
R)$ this implies
$$
|f||T_{k/b}f(t)||G_{k}(t)|\in L^{1}(\Bbb R ).
$$
Hence, by the Lebesgue Dominated Convergence Theorem,
$$
\lim_{\ell \rightarrow \infty}\frac{1}{b}\int_{\Bbb R}
\overline{f_{\ell}(t)}f_{\ell}(t-k/b) G_{k}(t)\ dt =
$$
$$
\frac{1}{b}\int_{\Bbb R}
\overline{f(t)}f(t-k/b) G_{k}(t)\ dt.
$$
Finally, since the right hand side of the WH-frame identity has only a
finite number of non-zero k's (and the same ones for $f$ and all  
$f_{\ell}$),
we have the equality in Step 3 and unconditional convergence in the
right hand side of the Identity.

$(2)\Rightarrow (1)$:  For any $f$ supported on an interval of  
length 1/b,
we are assuming the WH-frame Identity holds.  But, $F_{2}(f)=0$ for all
such $f$.  So
$$
F_{1}(f) = b^{-1}\int_{\Bbb R}|f(t)|^{2}\sum_{n}|g(t-na)|^{2}\ dt <  
\infty.
$$
This implies that $G_{0}$ is bounded.  To see this, let
$I=[c,d]$ be any interval of
length $< 1/b$.  It suffices to show that $G_{0}$ is bounded here since
$G_{0}$ is a-periodic means it is bounded if it is bounded on all
intervals of any fixed length.  Let
$$
A_{n} = \{t\in I : |G_{0}(t)|\le n \}.
$$
Let $T_{n}:L^{2}[c,d]\rightarrow L^{2}[c,d]$ be given by
$T_{n}f = {\chi}_{A_{n}}f\cdot \sqrt{G_{0}}$.  The $T_{n}$ are  
bounded linear
operators and the family is pointwise bounded by the above.  Hence
they are uniformly bounded and so
$$
Tf = f\cdot G_{0}
$$
is a bounded linear operator.  But the norm of this ``multiplication''
operator is ess sup $|G_{0}(t)|$.

\qed \vspace{14pt}

We remark that we could simplify the proof of Theorem
\ref{CF} if $g\in$PF.  For in this case we
use the frame operator $S$ to get some of the needed convergence.
For example, in this case we would observe:
$$
\sum_{m,n\in \mathbb Z}|<f,E_{mb}T_{na}g>|^{2} = <Sf,f>,
$$
while
$$
\lim_{{\ell}\rightarrow \infty}<Sf_{\ell},f_{\ell}> = <Sf,f>.
$$

\section{Types of Convergence of the WH-Frame Identity}\label{3}
\setcounter{equation}{0}

Here we will consider when the WH-frame Identity holds with a stipulation
on the type of convergence of the infinite series on the right hand side
of the identity.  These results are variations on results of Casazza,
Christensen and Janssen \cite{CCJ}.  To extend the results of \cite{CCJ}
we need a well known fact which is really the Polarization Identity  
for $H$.

\begin{Prop}\label{PI}
For any operator $T$ on a complex Hilbert space $H$ we have for all  
$x,y\in H$:
$$
4<Tx,y> = <T(x+y),x+y> - <T(x-y),x-y> +
$$
$$
i<T(x+iy),x+iy>-i<T(x-iy,x-iy>.
$$
\end{Prop}

\begin{Cor}
If $T$ is an operator on a complex Hilbert space $H$, then
$$
\|T\| \le 2\text{sup}\{|<Tf,f>|:\|f\|\le 1\}.
$$
\end{Cor}

We need some notation for checking the convergence of the series
in the WH-frame Identity.

\begin{Def}
Let $g\in L^{2}(R)$ satisfy $G_{0}(t)\le B$ a.e. For any $f\in
L^{2}(R)$,
and any $K,L\in Z$, we let
$$
S_{K,L}f(t) =  \sum_{k=-L}^{K}f(t-k/b)G_{k}(t),
$$
where as usual $G_{k}(t) = \sum_{n}g(t-na)\overline{g(t-na-k/b)}$.   
We also
let $S_{K} = S_{K,K}$.  If $M\subset Z$ with $|M|<\infty$, define
$$
S_{M}f(t) = \sum_{k\in M}f(t-k/b)G_{k}(t).
$$
\end{Def}

If $ \lim_{K \rightarrow \infty}  S_{K}f$ exists, we say that {\bf  
the Walnut
series for} $f$ {\bf
 converges symmetrically} - and this can be in either the norm or the
weak topology - and we say {\bf the Walnut series for f converges} when
$\lim_{K,L\rightarrow \infty}S_{K,L}f$ exists.

Now we give an extension of Theorem 5.2 from \cite{CCJ}.  The  
notation can
be found in Section \ref{1}.

\begin{Thm}[Casazza/Christensen/Janssen]\label{CCJS}
Let $a,b\in \mathbb R$, $g\in \L$ and $G_{0}(t)\le B$ a.e.  The following
are equivalent:

(1)  The Walnut series converges in norm symmetrically for every  
$f\in \L$.

(2)  The Walnut series converges weakly summetrically for every  
$f\in \L$.

(3)  We have $\text{sup}_{K}\|S_{K}\|<\infty$.
\\
Moreover, in this case the WH-system $(g,a,b)$ has a finite upper frame
bound and the Walnut series converges symmetrically to $Sf$, for all 
$f\in \L$.
\end{Thm}

We now extend this result slightly.

\begin{Thm}\label{CL1}
Let $a,b\in \mathbb R$, $g\in \L$ and $G_{0}(t)\le B$ a.e.  The following
are equivalent:

(1)  The Walnut series converges in norm symmetrically for every  
$f\in \L$.

(2)  We have for all $f\in \L$,
$$
\lim_{K\rightarrow \infty}<S_{K}f,f> = <Sf,f>,
$$

(3)  We have for all $f\in \L$,
$$
\sum_{m,n\in \mathbb Z}|<f,E_{mb}T_{na}g>|^{2} =
$$
$$
\lim_{K\rightarrow \infty}\frac{1}{b}\sum_{k=-K}^{K}\int_{\mathbb R}
\overline{f(t)}f(t-k/b)G_{k}(t)\ dt.
$$
\end{Thm}

{\it Proof.}
$(1)\Leftrightarrow (2)$:  If we assume (2), we easily obtain from  
Proposition
\ref{PI} that $\lim_{K\rightarrow \infty}<S_{K}f,h> = <Sf,h>$, for all
$f,h\in \L$.  i.e.  The Walnut series converges weakly symmetrically, and
hence symmetrically in norm.  The converse is obvious.

$(3)\Rightarrow (2)$:  The right hand side of (3) is:
$$
\lim_{K\rightarrow \infty}\frac{1}{b}\sum_{k=-K}^{K}\int_{\mathbb R}
\overline{f(t)}f(t-k/b)G_{k}(t)\ dt = \frac{1}{b}
\lim_{K\rightarrow \infty}<S_{K}f,f>.
$$
This implies that
$$
\lim_{K\rightarrow \infty}<S_{K}f,f> \ \ \text{exists}
$$
for all $f\in \L$.  By Proposition \ref{PI}, it follows that  
$(S_{K}f)$ is
weakly symmetrically
convergent (and hence symmetrically norm convergent by
Theorem \ref{CCJS}) for all $f\in \L$.  Hence,
the $(S_{K})$ are uniformly bounded operators.
By Proposition \ref{PP}, the right hand side of (3) converges to $<Sf,f>$
unconditionally
on a dense subset of $\L$, and since the operators $S_{K}$ are uniformly
bounded, we have the equality in (2) for all $f\in \L$.

$(1)\Rightarrow (3)$:  By (1), the limit on the right hand side of
(3) converges for all $f$ and to $<Sf,f>$.   By Theorem \ref{CCJS},  
$(g,a,b)$ has a finite
upper frame bound.  Now,
$$
\lim_{K\rightarrow \infty}\frac{1}{b}\sum_{k=-K}^{K}\int_{\mathbb R}
\overline{f(t)}f(t-k/b)G_{k}(t)\ dt = \lim_{K\rightarrow  
\infty}<S_{K}f,f> =
$$
$$
<Sf,f> = \sum_{m,n\in \mathbb Z}|<f,E_{mb}T_{na}g>|^{2}.
$$
\qed \vspace{14pt}

In \cite{CCJ} (Example 5.4) it is shown that there are WH-frames  
$(g,1,1)$ so
that for some function $f\in \L$, the Walnut series for $f$
does not converge symmetrically.  Combined with Theorem \ref{CL1}  
we obtain,

\begin{Cor}
There is a WH-frame $(g,1,1)$ and a function $f\in \L$ so that the
WH-frame identity fails for this $f$ in the sense that the series
on the right hand side of the WH-frame Identity does not converge
symmetrically for this $f$.
\end{Cor}

Next we generalize another result, Theorem 5.5, from \cite{CCJ}.

\begin{Thm}\label{CL2}
Let $ab\in \mathbb R$ and $g\in $PF.  The following are equivalent:

(1)  The Walnut series converges in norm for every $f\in \L$.

(2)  The Walnut series converges weakly for every $f\in \L$.

(3)  We have that $\text{sup}_{K,L}\|S_{K,L}\|< \infty.$

(4)  We have for all $f\in \L$,
$$
\lim_{K,L\rightarrow \infty}<S_{K,L}f,f> = <Sf,f>.
$$

(5)  We have for all $f\in \L$,
$$
\sum_{m,n\in \mathbb Z}|<f,E_{mb}T_{na}g>|^{2} =
$$
$$
\lim_{K,L\rightarrow \infty}\frac{1}{b}\sum_{k=-L}^{K}\int_{\mathbb R}
\overline{f(t)}f(t-k/b)G_{k}(t)\ dt.
$$
\end{Thm}

{\it Proof.}
The equivalence of $(1)-(3)$ is due to Casazza, Christensen, and Janssen
 (\cite{CCJ}, Theorem 5.5).  The rest of the proof follows line by  
line the
proof of our Theorem \ref{CL1} above, just replacing symmetric  
convergence by
convergence at each step.
\qed \vspace{14pt}

Again, in \cite{CCJ} (Example 5.7) it is shown that there is a WH-frame
$(g,1,1)$ for which the Walnut series converges symmetrically for every
$f\in \L$, but there is an $h\in \L$ for which the Walnut series does
not converge in norm (or weakly).  Combined with Theorem \ref{CL2},  
we have,

\begin{Cor}
There is a WH-frame $(g,1,1)$ for which the WH-frame Identity holds for
all $f\in \L$ in the sense that the series on the right hand side of the
identity converges symmetrically for all $f\in \L$ and we have  
equality in
the identity.  However, there is an $h\in \L$ for which the series on
the right hand side of the WH-frame Identity does not converge.
\end{Cor}

Our next theorem again generalizes  a result
(Theorem 6.1) from \cite{CCJ} and the proof follows
along the lines of the proof of Theorem \ref{CL1}.

\begin{Thm}\label{CL3}
Let $a,b\in \mathbb R$ and $g\in$PF.  The following are equivalent:

(1)  The Walnut series converges weakly unconditionally for every  
$f\in \L$.

(2)  The Walnut series converges unconditionally in norm for every  
$f\in \L$.

(3)  We have $\text{sup}_{M\subset \mathbb  
Z,|M|<\infty}\|S_{M}\|<\infty.$

(4)  We have that the series
$$
\sum_{k}<(T_{k/b}f)G_{k},f>
$$
converges unconditionally to $<Sf,f>$, for all $f\in \L$.

(5)  The WH-frame Identity holds and the series on the right hand side
converges unconditionally for all $f\in \L$.
\end{Thm}

In \cite{CCJ} (Example 6.3) it is shown that there is a WH-frame  
$(g,1,1)$
so that for every $f\in \L$ the Walnut series for $f$ converges in norm,
but there is some $h\in \L$ for which the Walnut series does not converge
unconditionally.  Combined with Theorem \ref{CL3} we have,

\begin{Cor}
There is a WH-frame $(g,1,1)$ so that for all $f\in \L$, the series on
the right hand side of the WH-frame Identity converges and is equal to
the left hand side.  However, there is a function $h\in \L$ so that the
series on the right hand side of the WH-frame Identity does not converge
unconditionally.
\end{Cor}

 Casazza, Christensen and Janssen \cite{CCJ} Theorem 6.5 have shown that
if $(g,a,b)$ satisfies the CC-condition, then for all $f\in \L$, the
Walnut series converges unconditionally.  Also, it is immediate that
if $g\in W(L^{\infty},L^{1})$, then $(g,a,b)$ satisfies the CC-condition
for all $a,b\in \mathbb R$.  These results, combined
with Theorem \ref{CL3} yields,

\begin{Cor}
If $(g,a,b)$ satisfies the CC-condition, then the WH-frame Identity holds
for all $f\in \L$ and the series converges unconditionally.  In  
particular,
if $g\in W(L^{\infty},L^{1})$ then the WH-frame Identity holds
for all $f\in \L$ and the series converges unconditionally.
\end{Cor}

\end{document}